\documentclass[oneside]{amsart}
\usepackage{amsmath,amsfonts,amssymb,amscd,amsbsy,amsthm,array}
\usepackage[mathscr]{eucal}
\usepackage{relsize}
\usepackage[plainpages=false,pdfpagelabels]{hyperref}
\usepackage[T1]{fontenc}
\usepackage{xspace,float}
\usepackage{rotating}

\DeclareMathOperator{\End}{End}

\newcommand{\sym}[1]{{\mathscr{S}_{#1}}}

\newcommand{\symm}[1]{{\sym{#1}}}

\newcommand{\coloneqq}{:\hspace{0pt}=}

\newcommand{\partition}{\vdash}

\newcommand{\ud}{\mathrm{d}}

\newcommand{\Complex}{\mathbb{C}}

\newcommand{\Integers}{\mathbb{Z}}

\newcommand{\QQ}{\mathbb{Q}}

\newcommand{\myexp}[1]{\exp\left({#1}\right)}

\newcommand{\mink}{{\textup{min}_k}}
\newcommand{\maxk}{{\textup{max}_k}}

\newcommand{\myDelta}[1]{\Delta \left( #1 \right)}
\newcommand{\myDeltaEmpty}{\Delta}

\newcommand{\myFatDelta}[1]{\myFatDeltaEmpty \left( #1 \right)}
\newcommand{\myFatDeltaEmpty}{\pmb{\Delta}}

\newcommand{\trivial}{\varnothing}

 \newtheorem{thm}{Theorem}
 \newtheorem{lem}[thm]{Lemma}
 \newtheorem{prop}[thm]{Proposition}
 \newtheorem{cor}[thm]{Corollary}

\newcommand{\FatP}{\mathbf{p}\vphantom{\tilde{\mathbf{p}}}}
\newcommand{\FatH}{\mathbf{h}\vphantom{\tilde{\mathbf{h}}}}
\newcommand{\FatSchur}{\mathbf{s}\vphantom{\tilde{\mathbf{s}}}}

\newcommand{\FatTildeP}{\tilde{\mathbf{p}}}
\newcommand{\FatTildeH}{\tilde{\mathbf{h}}}
\newcommand{\FatTildeSchur}{\tilde{\mathbf{s}}}

\newcommand{\TildeH}{\tilde{h}}
\newcommand{\TildeP}{\tilde{p}}

\newcommand{\DelP}[1]{\partial_{p_{#1}}}
\newcommand{\TildeDelP}[1]{\partial_{\TildeP_{#1}}}

\newcommand{\FatDelP}[1]{\partial_{\FatP_{#1}}}
\newcommand{\FatTildeDelP}[1]{\partial_{\FatTildeP_{#1}}}

\newcommand{\setOfX}{{\overline{X}}}
\newcommand{\setOfY}{{\overline{Y}}}

\newcommand{\myeval}{E_\sigma}
\newcommand{\myforget}{F}
\newcommand{\FatLambda}{\pmb{\Lambda}}

\newcommand{\BD}[2]{\textup{BD}^{#1 \, #2}}
\newcommand{\TW}[2]{\textup{TW}^{#1 \, #2}}
\newcommand{\FatBD}[2]{\textup{\textbf{BD}}^{#1 \, #2}}
\newcommand{\FatTW}[2]{\textup{\textbf{TW}}^{#1 \, #2}}
\newcommand{\Ratio}[2]{\textup{R}^{#1 \, #2}}

\newlength{\wdth}

\newcommand{\TWhT}[2]{\text{\makebox[20pt]{$\FatTildeH_{#2\mu_{#1}}$}}}
\newcommand{\TWhTnoadjust}[2]{\FatTildeH_{#2\mu_{#1}}}
\newcommand{\TWh}[2]{\text{\makebox[20pt]{$\FatH_{#2\lambda_{#1}}$}}}
\newcommand{\TWhnoadjust}[2]{\FatH_{#2\lambda_{#1}}}

\newcommand{\Raw}{\succ}
\newcommand{\Daw}{\curlyvee}
\newcommand{\Law}{\prec}

\newcommand{\mynwarrow}{\begin{turn}{135}{$ \leadsto $}\end{turn}}
\newcommand{\mysearrow}{\begin{turn}{-45}{$ \leadsto $}\end{turn}}

\begin{document}

\begin{abstract}
A classical question for a Toeplitz matrix with given symbol is to compute asymptotics for the determinants of its reductions to finite rank. One can also consider how those asymptotics are affected when shifting an initial set of rows and columns (or, equivalently, asymptotics of their minors). Bump and Diaconis~(\emph{Toeplitz minors}, J.~Combin.~Theory~Ser.~A, 97 (2002), pp.~252--271) 
 obtained a formula for such shifts involving Laguerre polynomials and sums over symmetric groups. They also showed how the Heine identity extends for such minors, which makes this question relevant to Random Matrix Theory. Independently, Tracy and Widom~(\emph{On the limit of some Toeplitz-like determinants}, SIAM~J.~Matrix~Anal.~Appl., 23 (2002), pp.~1194--1196)  used the Wiener-Hopf factorization to express those shifts in terms of products of infinite matrices. We show directly why those two expressions are equal and uncover some structure in both formulas that was unknown to their authors. We introduce a mysterious differential operator on symmetric functions that is very similar to vertex operators. We show that the Bump-Diaconis-Tracy-Widom identity is a differentiated version of the  classical Jacobi-Trudi identity.
\end{abstract}

\author{Paul-Olivier Dehaye}
\address{Merton College, University of Oxford, United Kingdom}
\email{paul-olivier.dehaye@maths.ox.ac.uk}
\date{\today}

\title[{On an identity due to ((Bump and Diaconis) and (Tracy and Widom))}]{On an identity due to \\{\large(}(Bump and Diaconis) and (Tracy and Widom){\large)}}

\maketitle

\section{Introduction}
\subsection{Origin: Toeplitz determinants}
Fix $\sigma(t)$ to be a function of the unit circle $\mathbb{T}$ in
$\Complex$ that can be written in the form
$$\sigma(t)=\exp \left(\sum_{k>0} \frac{p_k}{k} t^k +
\frac{\TildeP_k}{k}t^{-k}\right)$$
 for the sets of constants $\{p_k \in \Complex\}$ and  $\{\TildeP_k \in
\Complex \}$\footnote{This implies that $\sigma(0)=1$, a condition that is merely there for exposition.}. This requires $\sigma$ to have winding number~0 around 
the origin (since $\log \sigma(t)$ is defined, see~\cite[pp. 15--17]{BS} for more details). This also defines a set of constants $\{ d_k \}$ so that
$\sum_{k \in \Integers} d_k t^k \coloneqq 
\sigma(t)$ (i.e. the $d_k$'s are the Fourier coefficients of $\sigma(t)$). 
We will further assume that the $|p_k|$'s and $|\TildeP_k|$'s decrease fast enough, i.e. that all of the sums $\sum_k \frac{|p_k|}{k}, \sum_k \frac{|\TildeP_k|}{k}$ and $\sum_k \frac{|p_k \TildeP_k|}{k}$ are bounded.

We now construct a matrix $M_n$ having constant entries on diagonals
parallel to the main diagonal 
(Toeplitz property with \emph{symbol} $\sigma$):
$$
M_n(\sigma) = M_n = 
\begin{pmatrix}
d_0 & d_1 & \cdots & \cdots & d_{n-1}\\
d_{-1} & d_0 & d_1 & \cdots & d_{n-2}\\
\vdots & \ddots & \ddots & \ddots &\vdots\\
\vdots & \ddots & \ddots & \ddots &d_1\\
d_{1-n}& \cdots & \cdots & d_{-1} & d_0
\end{pmatrix}
_{n\times n} = \left( d_{i-j} \right)_{n\times n}
$$

A classical question for Toeplitz matrices is then to consider the
asymptotics of the determinant $\det(M_n)$ as $n$ goes to
infinity. Our identity will stem from the same question for a slightly
altered version of $M_n$. 

For $\lambda$ and $\mu$ partitions of length less or equal to $n$, look at 
$$ M_n^{\lambda \, \mu}(\sigma) := \left( d_{\lambda_i-\mu_j-i+j} \right)_{n \times n}.$$
Those new matrices are not Toeplitz, but at least they are minors of the Toeplitz matrix $M_m(\sigma)$, for some $m$ larger than $n$.  This is  clear once illustrated. For the sake of example, set $n:= 3$, $m:=5$,  $\lambda := (2,1)$, $\mu := (1)$. We then have the matrices
$$
M_3^{\lambda \, \mu}(\sigma) = \begin{pmatrix}
d_1&d_3&d_4\\
d_{-1} & d_1 & d_2 \\
d_{-3}&d_{-1}&d_0\\ 
\end{pmatrix}
\quad\text{and}\quad
M_5(\sigma) = 
\begin{pmatrix}
d_0&d_1&d_2&d_3&d_4\\
d_{-1}&d_{0}&d_{1}&d_{2}&d_{3}\\
d_{-2}&d_{-1}&d_{0}&d_{1}&d_{2}\\
d_{-3}&d_{-2}&d_{-1}&d_{0}&d_{1}\\
d_{-4}&d_{-3}&d_{-2}&d_{-1}&d_{0}
\end{pmatrix}.
$$
Observe that $M_3^{\lambda \, \mu}(\sigma)$ is the minor of $M_5(\sigma)$ obtained by striking its first and third column and its second and fourth row. If $m$ had been bigger, we would only have needed to strike more rows and columns.

The asymptotics of the determinants of $M^{\lambda \mu}(\sigma)$ are well known through the Szeg\"o limit theorem, so it is natural to look at the ratios
$$
R^{\lambda\mu}(\sigma):= \lim_{n \rightarrow \infty} \frac{\det M_n^{\lambda \, \mu}(\sigma)}{\det M_n(\sigma)}.
$$
These ratios have indeed been studied by two pairs of researchers, independently.

Tracy and Widom~\cite{TW} obtained the asymptotics $R^{\lambda\mu}(\sigma)$ as determinants involving the Fourier coefficients in the Wiener-Hopf factorization 
\begin{eqnarray}
\sigma(t) &=& \exp\left(\sum_{k>0} \frac{p_k}{k}t^k\right) \cdot \exp\left(\sum_{k>0} \frac{\TildeP_k}{k}t^{-k}\right) \notag\\
&=:& \sum_{k \ge 0 } h_k t^k \cdot \sum_{k \ge 0 } \tilde{h}_k t^{-k}
\label{Wiener-Hopf}
\end{eqnarray}
of $\sigma(t)$. The second line serves as definition of the $h_k$s and $\tilde{h}_k$s. We will present the full expression they obtain in Equation~(\ref{defTW}). Meanwhile, we refer to that expression as $\TW{\lambda}{\mu}(\sigma)$. 

Bump and Diaconis~\cite{BD} generalized instead the Heine identity. This classical identity gives 
$$
\det M_n(\sigma) \sim_{n\rightarrow \infty} \int_{U(n)} \sigma(g) \, \ud g,
$$
with $\sigma(g) := \prod \sigma (t_i)$, $t_i$ being the eigenvalues of $g$.  They extended this to
 $$\det M_n^{\lambda \,\mu}(\sigma) \sim_{n\rightarrow \infty} \int_{U(n)} \sigma(g) s_\lambda(g) \overline{s_\mu(g)}\, \ud g,$$
with $s_\lambda,s_\mu$ the usual Schur polynomials applied to the eigenvalues of $g$. Thus all results presented here for Toeplitz matrices apply for twisted integrals as well (hence the interest for Random Matrix Theory), and Bump and Diaconis derived independently from Tracy and Widom a second expression  (presented in Section~\ref{BDsection}) for the following limit:
$$
 \BD{\lambda}{\mu}(\sigma) := 
  \lim_{n \rightarrow \infty} \frac{\int_{U(n)} \sigma(g) s_\lambda(g) \overline{s_\mu(g)}\, \ud g}{\int_{U(n)} \sigma(g)\, \ud g}.
$$

Tracy and Widom's theorems are valid under slightly more general conditions than Bump and Diaconis'. Lyons \cite{SzegoLyons} discusses this point  in detail.

We now wish to state the theorem alluded to in the title of this article.
\begin{thm}[\cite{BD,TW}]
\label{main}
Let $\lambda, \mu$ be partitions. Then, for $\sigma(t)$ such that $\sum_k \frac{|p_k|}{k},$ $\sum_k \frac{|\TildeP_k|}{k}$ and $\sum_k \frac{|p_k
\TildeP_k|}{k}$ are bounded, we have
$$\BD{\lambda}{\mu}(\sigma) = \Ratio{\lambda}{\mu}(\sigma) =
\TW{\lambda}{\mu}(\sigma).$$
\end{thm}
The proof of this theorem thus comes from two entirely disjoint papers.

\subsection{Concept}
Theorem~\ref{main} raises an immediate question. If one forgets its origins, Theorem~\ref{main} is a mysterious combinatorial identity $\BD{\lambda}{\mu}(\sigma) =  \TW{\lambda}{\mu}(\sigma)$. Our main goal for this paper will be to prove this identity more directly, without relying on Toeplitz determinants (i.e. $\Ratio{\lambda}{\mu}(\sigma)$). 

Let $\trivial$ be the trivial partition. We will show how this identity is a differentiated version of the Jacobi-Trudi identity. We proceed along the following stages:
\begin{enumerate}
\item Both $\BD{\lambda}{\mu}$ and $\TW{\lambda}{\mu}$ are functions of $\sigma$, but can also be seen as functions of the Fourier coefficients $\{ p_1,p_2,\cdots,\TildeP_1,\TildeP_2,\cdots \}$. Those functions turn out to be power series in those Fourier coefficients. This is present in \cite{BD} and partly in \cite{TW}.
\item As explained in Section~\ref{defs}, the variable set $\{ p_1,p_2,\cdots,\TildeP_1,\TildeP_2,\cdots \}$ can be replaced by $\{ \FatP_1, \FatP_2,\cdots,\FatTildeP_1,\FatTildeP_2,\cdots\}$, the union of the two sets of symmetric power sums in two separate sets of variables (say $\setOfX$ and $\setOfY$). Notationally, this will replace $\BD{\lambda}{\mu}$ and $\TW{\lambda}{\mu}$ with $\FatBD{\lambda}{\mu}$ and $\FatTW{\lambda}{\mu}$.
\item There are two related differential operators, $\myDeltaEmpty$ and $\myFatDeltaEmpty$, that act respectively on $\BD{\lambda}{\mu}$ or $\TW{\lambda}{\mu}$ and on $\FatBD{\lambda}{\mu}$ or $\FatTW{\lambda}{\mu}$ (see Section~\ref{secDiff}).
\item \begin{thm}
\label{BDderiv}
$$
 \myDelta{\BD{\lambda}{\trivial}\cdot
\BD{\trivial}{\mu}}=\BD{\lambda}{\mu} \quad \text{ and } \quad \myFatDelta{\FatBD{\lambda}{\trivial}\cdot
\FatBD{\trivial}{\mu}}=\FatBD{\lambda}{\mu}.
$$
\end{thm}

\item \begin{thm}
\label{TWderiv}
$$
    \myDelta{\TW{\lambda}{\trivial}\cdot 
\TW{\trivial}{\mu}}=\TW{\lambda}{\mu}\quad  \text{ and } \quad
    \myFatDelta{\FatTW{\lambda}{\trivial}\cdot 
\FatTW{\trivial}{\mu}}=\FatTW{\lambda}{\mu}.
$$\end{thm}

\item $\FatBD{\lambda}{\trivial} = \FatTW{\lambda}{\trivial} $ (Jacobi-Trudi identity, a classical identity in symmetric function theory).
\end{enumerate}
As we can see, everything is proved through analogues in symmetric function theory which specialize to the objects of original interest. This can only work by ignoring the Toeplitz determinant origin of the expressions $\BD{\lambda}{\mu}$ and $\TW{\lambda}{\mu}$, but still gives a (new) corollary about the structure of the determinants: 
\begin{cor}
\label{Rderiv}
$$
    \myDelta{\Ratio{\lambda}{\trivial}\cdot
\Ratio{\trivial}{\mu}}=\Ratio{\lambda}{\mu}.
$$
\end{cor}

\subsection{Organization of this paper}
In Section~\ref{defs}, we will review the notions of symmetric function
theory which we need.
In Section~\ref{BDsection}, we will define $\BD{\lambda}{\mu}$ and prove Theorem~\ref{BDderiv}. The next Section accomplishes the same for the Tracy-Widom side and Theorem~\ref{TWderiv}.
Section~\ref{lastStep} will be devoted to the proof of Theorem~\ref{main}.
We give in Section~\ref{lastSection} a couple of noteworthy
relations on the $\Ratio{\lambda}{\mu}$'s. Finally, we discuss in Section~\ref{conclusion} how this paper fits into a more general program.

This research is part of the author's Ph.D. thesis \cite{DehayeThesis} at Stanford
University. It was supported in part by NSF grant FRG DMS-0354662. The
author wishes to acknowledge his adviser, Prof.~Daniel Bump, and Prof.~Persi Diaconis for extended discussions, as well as Prof.~Bertfried Fauser for pointing out relevance to his work with Prof.~Peter Jarvis. The author also thanks Ashkan Nikeghbali and Julie Rowlett for encouragements.

\section{General definitions and notations}
\label{defs}
We summarize the definitions and notations employed in this paper.
\subsection{Partitions and symmetric groups}
A partition $\lambda=(\lambda_1, \lambda_2, \cdots, \lambda_n)$ is a finite
decreasing sequence of non-negative integers. We define the weight
$|\lambda|$ of $\lambda$ to be the sum $\sum \lambda_i$. 
If this weight is $k$, we also use the notation $\lambda \partition k$. If
$k=0$, we denote the trivial partition $(0,0,0,0,\cdots)$ by
$\trivial$. The length $l(\lambda)$ of $\lambda$ is the maximal $i$ such
that $\lambda_i \neq 0$. 

There is a partial ordering on partitions: $\lambda \subseteq \mu$ iff
$\lambda_i\le \mu_i$ for all~$i$. In a probable break of standard
notation, $\lambda(i)$ counts the number of $\lambda_j$'s equal to $i$, so that
$(i^{\lambda(i)})=(\lambda_1,\lambda_2,\cdots,\lambda_n)$. In an
even greater offense, if $\pi$ is a permutation, we will use $\pi(i)$ for
the number of elements of $i$'s in the cycle type of $\pi$, \textbf{not}
for the image of point $i$ under $\pi$ (with no risk of notational confusion in the whole paper).

As usual, partitions of fixed weight $k$ index conjugacy classes in the
symmetric group on $k$ points $\symm{k}$. We set $z_\lambda := \prod_i
i^{\lambda(i)} i!$. This is the order of the centralizer of a permutation
in $\symm{|\lambda|}$ of cycle-type $\lambda$. 

In order to present the formula of Bump and Diaconis, we will also need
the irreducible characters of the symmetric groups. For a fixed~$k$, all
irreducible representations of $\symm{k}$ are indexed by partitions of
weight~$k$ (see the book by Sagan \cite{Sagan} for a friendly
introduction). If $\lambda \partition k$, we will use $\chi^\lambda$ for
the character of the representation corresponding to~$\lambda$. 

\subsection{Symmetric functions} 
We now introduce a few functions in the graded algebras
$\FatLambda(\setOfX)$ and  $\FatLambda(\setOfY)$ of symmetric functions in countably many
independent variables $\setOfX:= \{ x_1,x_2,x_3,\cdots \}$ and $\setOfY:= \{ y_1,y_2,y_3,\cdots \}$ over $\mathbb{Q}$. 
The former can be most directly thought of as the ring of formal sums
$S(x_1,\cdots)$ of monomials in the variables $x_i$ that have the symmetry
property $S(x_{\rho(1)},x_{\rho(2)},\cdots)=S(x_1,x_2,\cdots)$ for all
$\rho \in \symm{\infty}$. The most classic reference on the topic is
Macdonald's book \cite[Sections 1.2-1.5]{Macdonald}.

We will use the notation $\FatP_\lambda$,
$\FatH_\lambda$, $\FatSchur_\lambda$ and
$\FatSchur_{\lambda/\mu}$ for the various interesting functions
living in $\FatLambda(\setOfX)$. They will be respectively the power sum,
complete, Schur and skew Schur functions in the variables $\{x_i\}$
associated to the partition $\lambda$ (to the skew partition $\lambda/\mu$
for the latter). Similarly, we use $\FatTildeP_\lambda$,
$\FatTildeH_\lambda$, $\FatTildeSchur_\lambda$ and
$\FatTildeSchur_{\lambda/\mu}$ for the same functions in
$\FatLambda(\setOfY)$. We remind the reader that  boldface font will be used for
functions in $\FatLambda(\cdot)$. A tilde indicates the variable set $\setOfY$, while the default (when there is no tilde) is to assume that the variable set is $\setOfX$. 

One can define an inner product on $\FatLambda(\cdot)$ by setting the Schur polynomials to be orthonormal: $\left< \FatSchur_\lambda , \FatSchur_\mu \right>_{\FatLambda(\setOfX)}  = \delta_{\lambda \mu}$. 
The $\left< \cdot , \cdot \right>_{\FatLambda(\setOfX)}$ indicates that this inner product is for $\FatLambda(\setOfX)$. 
We will need the fact that the $\FatP_\lambda$'s form an orthogonal base: $\left< \FatP_{\lambda},\FatP_{\mu}\right>_{\FatLambda(\setOfX)} = z_\lambda \delta_{\lambda \mu}$.

We will also need to consider the algebra of symmetric functions in two sets of variables
$$
\FatLambda(\setOfX,\setOfY) = \FatLambda(\setOfX) \otimes_{\mathbb{Q}} \FatLambda(\setOfY).
$$
This comes equipped with an induced inner product defined by extending linearly  
$$
\left< \mathbf{a} \cdot \mathbf{\tilde{a}} , \mathbf{b} \cdot \mathbf{\tilde{b}} \right>_{\FatLambda(\setOfX,\setOfY)} = 
\left< \mathbf{a} , \vphantom{\tilde{b}}\mathbf{b} \right>_{\FatLambda(\setOfX)}\cdot
\left<  \mathbf{\tilde{a}} , \mathbf{\tilde{b}} \right>_{\FatLambda(\setOfY)}.
$$

\subsection{The derivations $\textbf{p}_n^\perp$ and $\tilde{\textbf{{p}}}_n^\perp$}
Let us first consider just the set of variables~$\setOfX$. 

Following Macdonald \cite[Example 3, Section 1.5, page 75]{Macdonald}, we define the algebra homomorphism $^\perp : \FatLambda(\setOfX) \longrightarrow \End(\FatLambda(\setOfX))$ in such a way that 
$$
\left< \textbf{f}^\perp \textbf{u},\textbf{v} \right>_{\FatLambda(\setOfX)} = \left< \textbf{u}, \textbf{f}\textbf{v} \right>_{\FatLambda(\setOfX)}   
$$
for all $\textbf{u},\textbf{v} \in \FatLambda(\setOfX)$. This is the adjoint of multiplication in the algebra $\FatLambda(\setOfX)$. 

Macdonald (following Foulkes) shows that $\FatP_n^\perp = n {\partial_{\FatP_{n}}}$ and so that $\FatP_n^\perp$ is a derivation. Indeed, we have
\begin{eqnarray*}
\left< \FatP_n^\perp( \FatP_\lambda) , \FatP_\mu \right>_{\FatLambda(\setOfX)} &=& \left< \FatP_\lambda , \FatP_\mu \FatP_n \right>_{\FatLambda(\setOfX)} \\
&=& 
\left\{
\begin{array}{cl}
0 & \text{ if } \lambda \ne (n) \cup \mu\\
z_\lambda & \text{ if } \lambda = (n) \cup \mu\\
\end{array}
\right.\\
&=& 
\left\{
\begin{array}{cl}
0 & \text{ if } \mu \ne \lambda \setminus (n)\\
z_\lambda & \text{ if } \mu \ne \lambda \setminus (n)\\
\end{array}
\right.\\
&=&\left<z_\lambda z_\mu^{-1} \FatP_{\lambda \setminus (n)} , \FatP_\mu\right>_{\FatLambda(\setOfX)}.
\end{eqnarray*}
But $z_\lambda z_{\lambda\setminus (n)}^{-1} = n\lambda(n)$, so 
$\FatP_n^\perp 
(\FatP_\lambda) = n \partial_{\FatP_{n}} (\FatP_\lambda)$ and we get our claim  that $\FatP_n^\perp = n {\partial_{\FatP_{n}}}$.

A similar result is of course true for $\FatLambda(\setOfY)$ (for the adjoint with respect to the inner product $\left< \cdot , \cdot \right>_{\FatLambda(\setOfY)}$).

Observe that 
$$
( \mathbf{a} \cdot \mathbf{\tilde{a}})^\perp = \mathbf{a}^\perp \otimes \mathbf{\tilde{a}}^\perp,
$$
and so $^\perp$ is a homomorphism $\FatLambda(\setOfX,\setOfY) \longrightarrow \End(\FatLambda(\setOfX,\setOfY))$.  
\label{Foulkes}

\subsection{Specializing symmetric objects}
Let $P = \left\{ p_1,p_2,\cdots \right\}$ and $\TildeP = \left\{ \TildeP_1,\TildeP_2,\cdots\right\}$ be sets of variables. We define $V_P = \QQ[[P]], V_{\TildeP} = \QQ[[\tilde{P}]]$ and $V = \QQ[[\tilde{P} \cup P]]$.

Any $\sigma(t)= \exp\left(\sum_{k>0} \frac{p_k}{k}t^k+ \frac{\TildeP_k}{k}t^{-k}\right)$ induces an evaluation map $\myeval: V\longrightarrow \Complex$ obtained by replacing the variables in $V$ by the values of the Fourier coefficients of $\log \sigma$. This is of course only convergent on a subset of $V$, but we will limit ourselves to that subset. 

We define algebra homomorphisms 
\begin{eqnarray*}
\myforget_\setOfX : \FatLambda(\setOfX) &\longrightarrow &V_P \subset V \quad(\text{resp. for }\setOfY, \tilde{P})\\
 \FatP_k &\longmapsto&p_k\\
\quad \quad \myforget : \FatLambda(\setOfX, \setOfY) &\longrightarrow &V\\
 \FatP_k &\longmapsto&p_k\\
 \FatTildeP_k &\longmapsto&\TildeP_k.
\end{eqnarray*}
Clearly, $\myforget$ restricts to $\myforget_\setOfX$ and $\myforget_\setOfY$, and merely forgets that the range was a vector space of symmetric polynomials: the variables $\setOfX$ and $\setOfY$ are completely lost. 

For a given $\sigma(t)$, we observe that the generating function for the $h_k$ is the same as the generating function for the $\FatH_k$, i.e. compare Equation~(\ref{Wiener-Hopf}) with the generating function identity 
\begin{eqnarray*}
\exp\left(  \sum_{k>0} \frac{\FatP_k }{k}t^k\right) = \sum_{k\ge0}
\FatH_k t^k.
\end{eqnarray*}

This very classical identity (Newton's identity describing the roots of a polynomial) was already discussed in the context of P\'olya's enumeration theory in the paper by Bump and Diaconis \cite{BD}. In any case, this guarantees that $$
\myeval \circ \myforget_\setOfX (\FatH_k) = h_k.
$$

Of course, a similar map $\myeval \circ \myforget_\setOfY: \FatLambda(\setOfY)\longrightarrow \Complex$ exists, and both maps together induce a third one, $\myeval \circ \myforget: \FatLambda(\setOfX,\setOfY)\longrightarrow \Complex$. 
We will call \emph{specialization} this whole process: start with a series in symmetric functions of countably many variables, forget through $\myforget$ that each symmetric function is a function itself (and thus assign a new variable for each function), and finally replace each of these new variables by a complex number through $\myeval$.

The advantage in setting up specialization in this way is that derivations are sent to derivations by $\myforget$. Thanks to Section~\ref{Foulkes}, we indeed know that for $\mathbf{f} \in \FatLambda(\setOfX,\setOfY)$, 
$$k \partial_{p_k} (\myforget (\mathbf{f})) = \myforget(k\partial_{\FatP_k}(\mathbf{f})) =\myforget(\FatP_k^\perp(\mathbf{f})).$$ 
We can use this property to create differential operators and specialize them from one algebra to another.
\subsection{\texorpdfstring{Differential operators  $\Delta$ and $\pmb{\Delta}$}{Differential operators}}  
\label{secDiff}
Consider still $V = \QQ[[P \cup \TildeP]]$. We define a (generalized) differential operator $\myDeltaEmpty $ as
$$
\myDeltaEmpty = \exp\left( \sum_k k \DelP{k} \TildeDelP{k} \right) =
\prod_{k > 0} \sum_{i \ge 0} \frac{k^i}{i!} (\DelP{k}\TildeDelP{k})^i,
$$
where $(\DelP{k}\TildeDelP{k})^i$ is composition. Note that sums and
products will be finite for any element of $V$ but that the order of
$\myDeltaEmpty$ is not uniformly bounded on~$V$. 

We define the operator $\myFatDeltaEmpty$ on $\FatLambda(\setOfX,\setOfY)$ in the same way (merely replacing $ \DelP{k} $ by $ \FatDelP{k} $). This implies the commutation relation 
\begin{eqnarray}
\myforget \circ \myFatDeltaEmpty &=&\myDeltaEmpty \circ \myforget.\label{eqnCommutation}
\end{eqnarray}

It follows from the previous sections that 
$$
\myFatDeltaEmpty =  \exp\left( \sum_k \frac{ \FatP_{k}^\perp \FatTildeP_{k}^\perp}{k} \right).
$$

\section{The Bump-Diaconis side}
\label{BDsection}
Assume $\lambda \partition m$ and $\mu \partition p$. Then Bump and
Diaconis define for each $\sigma(t) = \exp\left(\sum_{k>0} \frac{p_k}{k}t^k\right) \cdot \exp\left(\sum_{k>0} \frac{\TildeP_k}{k}t^{-k}\right)$ an expression $\BD{\lambda}{\mu}(\sigma) = \myeval(\BD{\lambda}{\mu})$, where they have
$$
\BD{\lambda}{\mu} = \frac{1}{m!}\sum_{\pi \in
S_m}\frac{1}{p!}\sum_{\rho \in S_p}
\chi^\lambda(\pi)\chi^\mu(\rho)\prod_{k>0}\text{F}_k(\pi,\rho),
$$
with 
$$
\text{F}_k(\pi,\rho)= 
\left\{
\begin{array}{cc}
k^{\rho(k)}\rho(k)!L^{(\pi(k)-\rho(k))}_{\rho(k)}\left( -\frac{p_k
\TildeP_k}{k}\right)p_k^{\pi(k)-\rho(k)} & \text{ if } \pi(k)\ge
\rho(k)\footnotemark\\
k^{\pi(k)}\pi(k)!L^{(\rho(k)-\pi(k))}_{\pi(k)}\left( -\frac{p_k
\TildeP_k}{k}\right)\TildeP_k^{\rho(k)-\pi(k)} & \text{ if } \rho(k)\ge
\pi(k)\\ 
\end{array}
\right.
\footnotetext{We remind the
reader of our unconventional usage of $\pi(k)$ for the number of $k-$cycles
in $\pi$.}
$$
and where
$$
L_n^{(\alpha)}(t) = \sum_{k=0}^n \binom{n+\alpha}{n-k} \frac{(-t)^k}{k!} =
\sum_{k=0}^n \binom{n+\alpha}{k}\frac{(-t)^{n-k}}{(n-k)!}
$$
is the usual Laguerre polynomial (the former expression is the standard 
definition, while the latter formula is only a reindexing of it that will be more useful here).

We define similarly $\FatBD{\lambda}{\mu}$ and $\mathbf{F}_k(\pi,\rho)$.

\begin{lem}
\label{betterBD}
Let $\maxk = \max(\pi(k),\rho(k))$ and  $\mink = \min(\pi(k),\rho(k))$. Then, 
$$
\mathbf{F}_k(\pi,\rho) = 
\sum_{i=0}^{\mink} k^i i! 
\binom{\maxk}{i}\binom{\mink}{i}
\FatP_k^{\pi (k)-i} \FatTildeP_k^{\rho(k)-i}.
$$
\end{lem}
\begin{proof}
We just need to expand the Laguerre polynomial in the definition of $\mathbf{F}_k$
while keeping track of the degrees in $\FatP_k$ and $\FatTildeP_k$. The key is to
observe that all the monomials will have the correct degrees, i.e. will be 
$\FatP_k^{\pi(k)-i}\FatTildeP_k^{\rho(k)-i}$ for $0\le i \le
\min(\rho(k),\pi(k))$. 
\end{proof}

\begin{proof}[Proof of Theorem \ref{BDderiv}]
When one of the partitions is trivial, the $\FatBD{\lambda}{\mu}$ reduce\footnote{We use here permutations as index for the power sums functions. We mean by $\FatP_\pi$ the function $\FatP_\lambda$, where $\lambda$ is the cycle-type of $\pi$.} to 
$$
\FatBD{\lambda}{\trivial}= \frac{1}{m!}\sum_{\pi \in \symm{m}} \chi^\lambda(\pi)\FatP_\pi
\quad \text{and} \quad \FatBD{\trivial}{\mu}=
\frac{1}{p!}\sum_{\rho \in \symm{p}} \chi^\mu(\rho)\FatTildeP_\rho
$$
We thus need to evaluate 
\begin{eqnarray*}
\myFatDelta{\FatBD{\lambda}{\trivial}\cdot \FatBD{\trivial}{\mu}} & = &
\myFatDelta{\frac{1}{m!}\sum_{\pi \in \symm{m}} \chi^\lambda(\pi)\FatP_\pi
\cdot
\frac{1}{p!}\sum_{\rho \in \symm{p}} \chi^\mu(\rho)\FatTildeP_\rho
}\\
&=&\frac{1}{m!}\sum_{\pi \in \symm{m}} \frac{1}{p!}\sum_{\rho \in
\symm{p}}  \chi^\lambda(\pi) \chi^\mu(\rho) \myFatDelta{\FatP_\pi \FatTildeP_\rho}.
\end{eqnarray*}
Each term is of the form
\begin{eqnarray*}
\myFatDelta{\FatP_\pi \FatTildeP_\rho} = \left[\prod_{k>0} e^{k\FatDelP{k}\FatTildeDelP{k}}\right]
\left(\prod_{k>0} \FatP_k^{\pi(k)}\FatTildeP_k^{\rho(k)}\right) = \prod_{k>0} \left[ e^{k\FatDelP{k}\FatTildeDelP{k}}
\left(\FatP_k^{\pi(k)}\FatTildeP_k^{\rho(k)}\right)\right], 
\end{eqnarray*}
where
\begin{eqnarray*}
\left[e^{k\FatDelP{k}\FatTildeDelP{k}}\right]
\left(\FatP_k^{\pi(k)}\FatTildeP_k^{\rho(k)}\right) & = & \sum_{i\ge 0}
\frac{(k\FatDelP{k}\FatTildeDelP{k})^i}{i!}
\left(\FatP_k^{\pi(k)}\FatTildeP_k^{\rho(k)}\right) \\
&=&\sum_{i \ge 0} k^i i! \binom{\pi(k)}{i} \FatP_k^{\pi(k)-i} \binom{\rho(k)}{i} \FatTildeP_k^{\rho(k)-i} \\
&=& \mathbf{F}_k(\pi,\rho)
\end{eqnarray*}
by Lemma~\ref{betterBD}.

Summing over all terms, we have
\begin{eqnarray*}
\myFatDelta{\FatBD{\lambda}{\trivial}\cdot \FatBD{\trivial}{\mu}} & = &
\frac{1}{m!}\sum_{\pi \in \symm{m}} \frac{1}{p!}\sum_{\rho \in \symm{p}}
\chi^\lambda(\pi) \chi^\mu(\rho) \myFatDelta{\FatP_\pi \FatTildeP_\rho}\\
& = & \frac{1}{m!}\sum_{\pi \in \symm{m}} \frac{1}{p!}\sum_{\rho \in
\symm{p}}  \chi^\lambda(\pi) \chi^\mu(\rho) \prod_{k>0}\mathbf{F}_k(\pi,\rho)\\
& = & \FatBD{\lambda}{\mu}.
\end{eqnarray*}
The identity
$$
\myDelta{\BD{\lambda}{\trivial}\cdot \BD{\trivial}{\mu}}= \BD{\lambda}{\mu}
$$
follows from
\begin{eqnarray*}
\BD{\lambda}{\mu} &= & \myforget(\FatBD{\lambda}{\mu})\\
&=& \myforget\left(\myFatDelta{\FatBD{\lambda}{\trivial}\cdot \FatBD{\trivial}{\mu}}\right)\\
&=& \myDelta{\myforget(\FatBD{\lambda}{\trivial}\cdot \FatBD{\trivial}{\mu})} \quad\quad\text{(Equation~(\ref{eqnCommutation}))}\\
&=&\myDelta{\BD{\lambda}{\trivial}\cdot \BD{\trivial}{\mu}},
\end{eqnarray*}
completing our proof of Theorem \ref{BDderiv}.
\end{proof}

We will need an additional lemma later.
\begin{lem}
\label{SchurBD}
$$
\FatBD{\lambda}{\trivial} = \FatSchur_\lambda\quad \text{and} \quad
\FatBD{\trivial}{\mu}= \FatTildeSchur_\mu.$$
\label{SchurExpansion}
\end{lem}
\begin{proof}
This is immediate from the definitions of $\FatBD{\lambda}{\trivial}$ and 
$\FatBD{\trivial}{\mu}$: we get the expansions\footnote{Again, we use here permutations as index for the power sums functions.} 
$$\frac{1}{|\lambda|!}\sum_{\pi \in \symm{|\lambda|}} \chi^\lambda(\pi) \FatP_\pi = \FatSchur_\lambda \quad \text{and} \quad \frac{1}{|\lambda|!}\sum_{\pi \in \symm{|\lambda|}} \chi^\lambda(\pi) \FatTildeP_\pi = \FatTildeSchur_\lambda $$ 
for Schur polynomials in terms of power sums, a fact that was already
presented by Bump and Diaconis in their paper. 
\end{proof}

\section{The Tracy-Widom side}
\label{TWsection}
Since $\sigma(t)=\exp \left(\sum_{k>0} \frac{p_k}{k} t^k +
\frac{\TildeP_k}{k} t^{-k}\right)$, it is reasonable to consider the
functions
\begin{eqnarray*}
\sigma^+(t)  := & \sum_{k \ge 0} h_k t^k :=& \exp \left(\sum_{k>0} \frac{p_k}{k} t^k\right) \\
\text{and }\sigma^-(t)  := & \sum_{k \ge 0} \TildeH_k t^{-k} :=  & \exp \left(\sum_{k>0} \frac{\TildeP_k}{k}
t^{-k}\right). \\
\end{eqnarray*}
It is a classical theorem from operator theory for Toeplitz matrices (see
B\"ottcher and Silbermann's book \cite[page 15]{BS}) that we then have
$$
\lim_{n \rightarrow \infty} (M_n(\sigma^+) \cdot M_n(\sigma^-))_{ij} =
\lim_{n \rightarrow \infty} M_n(\sigma)_{ij}.
$$
This is called the Wiener-Hopf factorization of the symbol $\sigma$.

Tracy and Widom use the Fourier coefficients $h_k$'s and $\TildeH_k$'s of $\sigma^+$ and $\sigma^-$ to formulate their result.

We are now ready to define $\FatTW{\lambda}{\mu}$ for the
partitions $\lambda \partition m$ and $\mu \partition p$. Let $d$ be an
integer large enough that $\lambda_{d+1} = \mu_{d+1} = 0$. Obviously,
$d=\max ( l(\lambda), l(\mu))$ would do, but $d$ could be taken larger
without affecting the result. Then we set
\begin{eqnarray}
\nonumber
\FatTW{\lambda}{\mu} &:=& \det \left(\left(
\FatTildeH_{i-j+\mu_{d-i+1}}\right)_{d \times \infty} \cdot \left(
\FatH_{j-i+\lambda_{d-j+1}}\right)_{\infty\times d}\right)
\\
&=& \det \left(
\left(
\begin{smallmatrix}
\TWhTnoadjust{d}{} & & & \Raw & \TWhTnoadjust{d}{1-d+ } &
\TWhTnoadjust{d}{-d+} & \Raw &
                & & \Raw &  \FatTildeH_{0} &  0  & \cdots &  &  & & 0&
\cdots\\
 & \mynwarrow & & \Raw & & & \Raw &
                 &  & \Raw & \FatTildeH_{0} & 0 & \cdots & & & & 0 &
\cdots\\
 & & \TWhT{d-i+1}{} & &  & &\Raw &
                &  & & \Raw & \FatTildeH_{0} & 0 & \cdots & & & 0&
\cdots\\
 & & & \mynwarrow & & & \Raw &
                & & & \Raw &  \FatTildeH_{0} & 0 & \cdots & & & 0 &
\cdots\\
\TWhTnoadjust{1}{d-1+} & & & \Raw & \TWhT{1}{} & \TWhT{1}{-1+} & \Raw &
                & & & & & \Raw & \FatTildeH_{0} &  0  & \cdots & 0&
\cdots\\
\end{smallmatrix}
\right)_{d\times \infty} \right.\nonumber
\\
&&
\hspace{.5in}\cdot\hspace{.5in}
\left.
\left(
\begin{smallmatrix}
\TWh{d}{} & & & & \TWh{1}{d-1+}\\
\Daw&\mynwarrow &&&\Daw\\
& & \TWh{d-i+1}{}&&&\\
\Daw & & & \mynwarrow&\Daw\\
\TWh{d}{1-d+} & & & & \TWh{1}{}\\
\TWh{d}{-d+}   & & & & \TWh{1}{-1+}\\
\Daw & \Daw & \Daw & \Daw & \Daw\\
\\
\Daw&\Daw\\
\FatH_0&\FatH_0&\Daw&\Daw\\
0&0&\FatH_0&\FatH_0\\
\vdots&\vdots&0&0&\Daw\\
&&\vdots&\vdots&\FatH_0\\
&&&&\vdots\\
\\
0&0&0&0&0\\
\vdots&\vdots&\vdots&\vdots&\vdots\\
\end{smallmatrix}
\right)_{\infty \times d}
\right). \label{defTW}
\end{eqnarray}
The structure of those matrices is important. We now attempt to describe
it in words.

We have here the determinant of a product of two ``half-strip'' matrices of
sizes
$d\times\infty$ and $\infty \times d$. The entries along the main diagonal
(marked by the arrows $\mynwarrow$) are all of the form
$\FatH_{\lambda_i}$
or $\FatTildeH_{\mu_i}$. The first matrix (resp. second) has a privileged
direction, $\Raw$ (resp. $\Daw$), in which the indices of $\FatH_\star$
(resp. $\FatTildeH_\star$) are
decreasing. This guarantees that the product is well-defined: each line on
the first column has only
finitely many non-zero entries\footnote{This is not important, but there is also a ``cascading effect'' among non-zero entries: in the first matrix for instance, the last non-zero entry on each row (i.e.  $ \FatTildeH_0$) has to be (weakly) to the right of any non-zero entry on the rows above.}. 

We define $\TW{\lambda}{\mu} := \myforget ( \FatTW{\lambda}{\mu})$, and indeed the expression $\TW{\lambda}{\mu}(\sigma) = \myeval (\TW{\lambda}{\mu})$ is what appears in \cite{TW}. We make the pedantic distinction here between $\TW{\lambda}{\mu}$ and $\TW{\lambda}{\mu}(\sigma)$ to highlight that the former is an element of $V$, i.e. a power series in the variable set $P \cup \tilde{P}$, and can thus be differentiated, unlike the latter which is only a complex number. 

The matrices involved in the definitions of $\FatTW{\lambda}{\mu} $ and $ \TW{\lambda}{\mu}$ are obviously very similar to the Jacobi-Trudi matrix. 
We remind the reader that the Jacobi-Trudi matrix of dimension $d\times d$
for the partition $\lambda$ ($d \ge l(\lambda)$) is the matrix
$$
\textbf{JT}_\lambda^d = 
\begin{pmatrix}
\TWh{1}{} & \Law & & \Law & \TWhnoadjust{1}{d-1+}\\
&\mysearrow& & \Law&\\
&&\TWh{i}{}&&\\
&\Law&&\mysearrow&\\
\TWhnoadjust{d}{1-d+}& \Law&&\Law&\TWhnoadjust{d}{}\\
\end{pmatrix}_{d\times d},
$$
where we respected the same conventions with arrows. We define $\tilde{\textbf{JT}}^d_\lambda$ in a totally analogous way (i.e. using $ \FatTildeH$'s). It is a central theorem of the theory of symmetric functions that $\det
(\textbf{JT}_\lambda^d)=\FatSchur_\lambda$ (see \cite[Theorem~37.1]{BumpLieGroups}) and is thus independent of $d$ (as long as $d\ge l(\lambda)$). Similarly, $\det
(\tilde{\textbf{JT}}_\lambda^d)=\FatTildeSchur_\lambda$.

We are now ready to comment on the result of Tracy and Widom a bit further. 
\begin{lem}
\label{SchurTW}
$$
\FatTW{\lambda}{\trivial} = \FatSchur_\lambda \quad
\text{and}\quad \FatTW{\trivial}{\mu} = \FatTildeSchur_\mu.
$$
\end{lem}
\begin{proof}
We will only do the case $\mu = \trivial$. Pick $d \ge l(\lambda)$. The left-hand side matrix in the definition of $\FatTW{\lambda}{\trivial}$ is then lower triangular, with $1$'s on the main
diagonal. Without affecting the final determinant, we can row-reduce this
matrix to 
$
\left( \delta_{ij} \right)_{d \times \infty}  
$,
with $\delta_{ij}$ the Kronecker delta.

Hence we easily compute 
\begin{eqnarray*} \FatTW{\lambda}{\trivial} & =
&\begin{pmatrix}
\TWh{d}{} & & & & \TWh{1}{d-1+}\\
\Daw&\mynwarrow &&&\Daw\\
& & \TWh{d-i+1}{}&&&\\
\Daw & & & \mynwarrow&\Daw\\
\TWhnoadjust{d}{1-d+} & & & & \TWh{1}{}\\
\end{pmatrix}_{d\times d}\\
&=&\det \left( \left(
\textbf{JT}^d_\lambda\right)_{d+1-j,d+1-i}\right)=\det
\textbf{JT}^d_\lambda = \FatSchur_\lambda.
\end{eqnarray*}
The key observation is thus that the $d\times d$ truncation of the
right-hand side matrix in the Tracy-Widom determinant is the
anti-transpose\footnote{The anti-transpose of a matrix is its transposed
along the main anti-diagonal.} of the Jacobi-Trudi matrix, and
that a determinant is not affected under anti-transposition. 
\end{proof}

We can now get started on the proof of Theorem~\ref{TWderiv}.

\begin{proof}[Proof of Theorem~\ref{TWderiv}]
We need to compute $\myFatDelta{\FatTW{\lambda}{\trivial} \cdot
\FatTW{\trivial}{\mu} }$. We have
$$\myFatDeltaEmpty = \myexp{\sum_k k \FatDelP{k} \FatTildeDelP{k}} = \myexp{\sum_k \frac{\FatP_{k} \FatTildeP_{k}}{k} }^\perp.$$ 
\label{HeineToBe}
The exponential can easily be expanded to obtain 
$$\myFatDeltaEmpty = \left(\sum_\nu \frac{1}{z_\nu} \FatP_\nu \FatTildeP_\nu\right)^\perp,$$ where the sum is over all partitions $\nu$. We now make use of the Cauchy identity 
$$\sum_\nu \frac{1}{z_\nu} \FatP_\nu \FatTildeP_\nu = \prod_{\substack{x_i \in \setOfX\\ y_j \in \setOfY}} \frac{1}{1-x_i y_j} = \sum_\nu \FatSchur_\nu \FatTildeSchur_\nu$$ and obtain our final expression:
$$\myFatDeltaEmpty = \left(\sum_\nu \frac{1}{z_\nu} \FatP_\nu \FatTildeP_\nu\right)^\perp = \left(\sum_\nu \FatSchur_\nu \FatTildeSchur_\nu\right)^\perp. $$

Coming back to our original computation, we just obtained
\begin{eqnarray}
\myFatDelta{\FatTW{\lambda}{\trivial} \cdot \FatTW{\trivial}{\mu} } = \sum_\nu \FatSchur_\nu^\perp (\FatSchur_\lambda) \FatTildeSchur_\nu^\perp (\FatTildeSchur_\mu). \label{goal}
\end{eqnarray}
 Observe that
\begin{eqnarray*}
\FatSchur_\nu^\perp (\FatSchur_\lambda) & = & \sum_\mu \left< \FatSchur_\nu^\perp (\FatSchur_\lambda) , \FatSchur_\mu \right> \FatSchur_\mu \\
&=&\sum_\mu \left< \FatSchur_\lambda , \FatSchur_\mu \cdot \FatSchur_\nu \right> \FatSchur_\mu\\
&=&\sum_\mu c^\lambda_{\mu \nu} \FatSchur_\mu = \FatSchur_{\lambda/\nu}.
\end{eqnarray*}
The last sum, which involves the Littlewood-Richardson coefficients, is precisely the definition of $ \FatSchur_{\lambda/\nu}$. 

Armed with this observation, we can thus rework Equation~(\ref{goal}) into
$$
\myFatDelta{\FatTW{\lambda}{\trivial} \cdot \FatTW{\trivial}{\mu}} = \sum_\nu \FatSchur_{\lambda/\nu} \FatTildeSchur_{\mu/\nu}.
$$

When $\nu$ runs through all partitions, the skew function $\FatSchur_{\lambda/\nu}$ runs through all $d\times d $ minors 
$
\left(
\FatH_{j-i-\nu_i+\lambda_{d-j+1}}\right)_{d \times d}
$
of the matrix 
$
\left(
\FatH_{j-i+\lambda_{d-j+1}}\right)_{\infty\times d}
$. Similarly, $\FatTildeSchur_{\mu/\nu}$ will run through the minors 
$
\left(
\FatTildeH_{i-j-\nu_j+\mu_{d-i+1}}\right)_{d \times d}
$
of 
$
\left(
\FatTildeH_{i-j+\mu_{d-i+1}}\right)_{d \times \infty}
$. Moreover, the minors obtained with $\FatSchur_\nu^\perp$ and $\FatTildeSchur_\nu^\perp$ are paired up just as in the Cauchy-Binet identity. 
Therefore, we obtain
\begin{eqnarray*}
\myFatDelta{\FatTW{\lambda}{\trivial} \cdot \FatTW{\trivial}{\mu} } &=&\det \left(\left(
\FatTildeH_{i-j+\mu_{d-i+1}}\right)_{d \times \infty} \cdot \left(
\FatH_{j-i+\lambda_{d-j+1}}\right)_{\infty\times d}\right)
\\
&=&\FatTW{\lambda}{\mu}
\end{eqnarray*}
and we are done. The proof for $\TW{\lambda}{\mu}$ simply follows from applying the homomorphism~$\myforget$.
\end{proof}

\section{The proof of Theorem~\ref{main}}
\label{lastStep}
\begin{proof}
We have from Lemmas~\ref{SchurExpansion} and~\ref{SchurTW} that 
$$
\FatBD{\lambda}{\trivial}=\FatSchur_\lambda = \FatTW{\lambda}{\trivial}
\quad \text{and} \quad
\FatBD{\trivial}{\mu}=\FatTildeSchur_\mu = \FatTW{\trivial}{\mu}.
$$
Tracing back to those lemmas, this is a direct consequence of the Jacobi-Trudi identity. 

The Theorem now follows. We have
$$
\begin{array}{rcll}
\FatBD{\lambda}{\mu} &=
&\myFatDelta{\FatBD{\lambda}{\trivial}\cdot\FatBD{\trivial}{\mu}} &\quad
\text{(Theorem~\ref{BDderiv})}  \\
&=&\myFatDelta{\FatTW{\lambda}{\trivial}\cdot\FatTW{\trivial}{\mu}}&\quad
\text{(Lemmas~\ref{SchurBD} and \ref{SchurTW})}\\
&=&\FatTW{\lambda}{\mu} &\quad \text{(Theorem~\ref{TWderiv})}.

\end{array} 
$$ 
\end{proof}

\section{\texorpdfstring{Some relations among $\Ratio{\lambda}{\mu}$'s}{Some relations among ratios}}
\label{lastSection}
We now consider $\Ratio{\lambda}{\mu}$ as an element of $V$, and immediately see that Corollary~\ref{Rderiv} is a consequence of the previous theorems. 
Two very natural properties of $\Ratio{\lambda}{\mu}$ also pop out of the presentation due to Tracy and Widom. The proofs rely only on basic properties of determinants and the Tracy-Widom expression $\FatTW{\lambda}{\mu}$, and their statement does not involve differential operators. Unlike Corollary~\ref{Rderiv}, they could thus be stated by evaluation at a specific $\sigma$.
\begin{prop}
Let $(r)$ and $(s)$ denote partitions with just one part each, of size $r\ge 1$ and $s \ge 1$ and let $\lambda, \mu$ be partitions, with $\max (l(\lambda),l(\mu))\le d$. Then,
\begin{eqnarray}
\label{prop_TW_1}
\Ratio{(r)}{(s)} = \Ratio{(r)}{\trivial}\cdot \Ratio{\trivial}{(s)} + \Ratio{(r-1)}{(s-1)}
\end{eqnarray}
and
\begin{eqnarray}
\Ratio{\lambda}{\mu} = \det \left(\Ratio{(\lambda_i+d-i)}{(\mu_j+d-j)}\right)_{1\le i,j \le d}.
\label{prop_TW_2}
\end{eqnarray}
\end{prop}
\begin{proof}
Both results follow from the same fact: 
\begin{eqnarray*}
\FatTW{(r)}{(s)} &=& \det \left(\left(
\FatTildeH_{1-j+s}\right)_{1 \times \infty} \cdot \left(
\FatH_{1-i+r}\right)_{\infty\times 1}\right)
\\
&=& \FatTildeH_s \FatH_r + \FatTildeH_{s-1} \FatH_{r-1} + \cdots\\
&=& \FatTildeH_s \FatH_r + \FatTW{(r-1)}{(s-1)}\\
&=& \FatTW{(r)}{\trivial} \FatTW{\trivial}{(s)} + \FatTW{(r-1)}{(s-1)},
\end{eqnarray*}
which proves Equation~(\ref{prop_TW_1}).

For Equation~(\ref{prop_TW_2}), we just need to observe that $\FatTW{\lambda}{\mu}$ is defined as the determinant of a matrix $\textbf{M}$ which itself is a product of two matrices. The coefficient on the $i^\text{th}$ row and the $j^\text{th}$ column of $\textbf{M}$ is  given by 
$$
\textbf{M}_{ij} = \sum_{k=0}^\infty \FatTildeH_{i-1-k+\mu_{d+1-i}} \FatH_{j-1-k+\lambda_{d+1-j}},
$$
where this sum is actually finite (because the terms eventually vanish).

By the reasoning for Equation~(\ref{prop_TW_1}), we actually know that
$$
\textbf{M}_{ij} = \FatTW{(j-1+\lambda_{d+1-j})}{(i-1+\mu_{d+1-i})}.
$$

Equation~(\ref{prop_TW_2}) then follows from the invariance of determinants under transposition and anti-transposition.
\end{proof}

 \section{Conclusion and speculation}
 \label{conclusion}
To summarize this paper, we have reproved the (Bump-Diaconis)-(Tracy-Widom) identity (and proved Corollary~\ref{Rderiv}) through specialization from the deeper symmetric function identity
$$
\myFatDelta{\FatSchur_{\lambda}\cdot\FatTildeSchur_{\mu}} = \sum_\nu \FatSchur_{\lambda/\nu} \FatTildeSchur_{\mu/\nu}.
$$

We feel that this more axiomatic approach to random matrix theory integrals through the theory of symmetric functions has a lot of potential. The backbone of symmetric function theory is common with much of the work of Fauser and Jarvis \cite{HopfApproach, NewRules, HopfLaboratory} on group branchings (which they sometimes specialize for perturbative quantum field theory). In particular, the operator $\Delta$ appears as a twisted product or ``Cliffordization'' in \cite{HopfLaboratory}, and results generalizing theorems~\ref{BDderiv} and~\ref{TWderiv} to other groups have been obtained in \cite{NewRules}. Our methods however seem to be much simpler, mostly because the differential operator $\myDeltaEmpty$ allows to encode the Nywell-Littlewood formula in a generating series form. We thus hope the techniques presented here will naturally expand to all classical compact Lie groups. Note that this generalization for expressions of the type $R^{\trivial \lambda}$ has already been achieved (independently) in \cite{DehayeWeyl}.
\bibliographystyle{alpha}
\bibliography{../../../references/main.bib}

\begin{thebibliography}{FJKW06}

\bibitem[BD02]{BD}
Daniel Bump and Persi Diaconis.
\newblock Toeplitz minors.
\newblock {\em J. Combin. Theory Ser. A}, 97(2):252--271, 2002.

\bibitem[BS99]{BS}
Albrecht B{\"o}ttcher and Bernd Silbermann.
\newblock {\em Introduction to large truncated {T}oeplitz matrices}.
\newblock Universitext. Springer-Verlag, New York, 1999.

\bibitem[Bum04]{BumpLieGroups}
Daniel Bump.
\newblock {\em Lie groups}, volume 225 of {\em Graduate Texts in Mathematics}.
\newblock Springer-Verlag, New York, 2004.

\bibitem[Deh06]{DehayeThesis}
Paul-Olivier Dehaye.
\newblock {\em Averages over compact Lie groups, twisted by Weyl characters and
  application to moments of derivatives of characteristic polynomials}.
\newblock PhD thesis, Stanford University, 2006.
\newblock \url{http://www.maths.ox.ac.uk/~pdehaye/papers/thesis/thesis.pdf}.

\bibitem[Deh07]{DehayeWeyl}
Paul-Olivier Dehaye.
\newblock Averages over classical {L}ie groups, twisted by characters.
\newblock {\em J. Combin. Theory Ser. A}, 114(7):1278--1292, 2007.

\bibitem[FJ04]{HopfLaboratory}
Bertfried Fauser and Peter~D. Jarvis.
\newblock A {H}opf laboratory for symmetric functions.
\newblock {\em J. Phys. A}, 37(5):1633--1663, 2004.

\bibitem[FJK]{HopfApproach}
Bertfried Fauser, Peter~D. Jarvis, and Ronald~C. King.
\newblock A {H}opf algebraic approach to the theory of group branchings.
\newblock \url{arXiv:math-ph/0508034}.

\bibitem[FJKW06]{NewRules}
Bertfried Fauser, Peter~D. Jarvis, Ronald~C. King, and Brian~G. Wybourne.
\newblock New branching rules induced by plethysm.
\newblock {\em J. Phys. A}, 39(11):2611--2655, 2006.

\bibitem[Lyo03]{SzegoLyons}
R.~Lyons.
\newblock Szeg{\H o} limit theorems.
\newblock {\em Geom. Funct. Anal.}, 13(3):574--590, 2003.

\bibitem[Mac95]{Macdonald}
I.~G. Macdonald.
\newblock {\em Symmetric functions and {H}all polynomials}.
\newblock Oxford Mathematical Monographs. The Clarendon Press Oxford University
  Press, New York, second edition, 1995.
\newblock With contributions by A. Zelevinsky, Oxford Science Publications.

\bibitem[Sag01]{Sagan}
Bruce~E. Sagan.
\newblock {\em The symmetric group}, volume 203 of {\em Graduate Texts in
  Mathematics}.
\newblock Springer-Verlag, New York, second edition, 2001.
\newblock Representations, combinatorial algorithms, and symmetric functions.

\bibitem[TW02]{TW}
Craig~A. Tracy and Harold Widom.
\newblock On the limit of some {T}oeplitz-like determinants.
\newblock {\em SIAM J. Matrix Anal. Appl.}, 23(4):1194--1196 (electronic),
  2002.

\end{thebibliography}
\end{document}